\documentclass[final,12pt]{article}
\usepackage{verbatim}
\usepackage{latexsym}
\usepackage{amsmath,amsthm,amsfonts}
\usepackage{eucal}
\usepackage{cases}

\usepackage{mathtext}
\usepackage[cp1251]{inputenc}
\usepackage[T2A]{fontenc}
\usepackage[dvips]{graphicx}
\usepackage{amsmath}
\usepackage{amssymb}
\usepackage{amsxtra}
\usepackage{latexsym}
\usepackage{ifthen}
\usepackage{color}

\newtheorem{theorem}{Theorem}
\newtheorem{definition}{Definition}[section]

\newtheorem{example}{Example}[section]
\newtheorem{corollary}{Corollary}[section]
\numberwithin{equation}{section}

\begin{document}

\title{Thin film flow dynamics on fiber nets}

          \author{Roman M. Taranets$^a$,  Marina Chugunova$^b$ \\
          $^a$Institute of Applied Mathematics and Mechanics of the NASU,\\
1, Dobrovol'skogo Str., 84100, Sloviansk, Ukraine,\\
  taranets\_r@yahoo.com \\
        $^b$Institute of Mathematical Sciences, Claremont Graduate University,\\
150 E. 10th Str., Claremont, California 91711, USA, \\
 marina.chugunova@cgu.edu }

\maketitle

\begin{abstract}
We analyze existence and qualitative behavior of non-negative weak solutions for fourth order degenerate
parabolic equations on graph domains with Kirchhoff's boundary conditions at the inner nodes and Neumann boundary conditions at the boundary nodes. The problem is originated from industrial constructions of spray coated meshes which are used in water collection and in oil-water separation processes. For a certain range of parameter values we prove  convergence toward a constant steady state that corresponds to the uniform distribution of coating on a fiber net.
          \end{abstract}



\section{Introduction}

Creation of fiber nets for water collection from the air is a new research area in industrial engineering.
Fog represents a large source of drinkable water, and considered to be one of possible solutions during droughts
in arid climates. Different plants and animals developed special textural and chemical features on their surfaces
to harvest this resource of water. \cite{optdesign} investigated the influence of the
surface wettability characteristics, length scale, and weave density on the fog-harvesting
capability of woven meshes.

 Bio-inspired fibers have been researched significantly so as to create a new type of meshes in fogging-collection projects. \cite{bioinsp}
 confirmed that the water-collecting ability of the spider web is the result of a unique fiber structure. They  fabricated different types of meshes
 to investigate the water collection behavior and the influence of geometry on the hanging-drops and used these data
  to evaluate the length of the three phase contact line at threshold conditions in conjunction with the maximal volume
 of a hanging drop at different modes. They also demonstrated that the geometrical structure of spider-inspired fiber
 induced much stronger water hanging ability than that of uniform fiber.

  On a rainy day a spider web cannot avoid collision with water droplets. The collision normally does not destroy the fiber
  network of the web, and droplet residue is collected on the fiber after impact. After a series of impacts, the web is covered
  by a number of water drops, which are even larger than normal rain drops (see \cite[Figure 3.1]{rain-drop}).
  Drop impact on spiders web is also encountered in engineering problems because the water harvesting performance crucially depends
   on the water retention on the fibers. On the contrary,  in air filtration systems using meshes, droplet retention must be suppressed
   because it causes clogging of fiber filters (see \cite{clogging}).

    When products in the form of a mesh or grid are coated by spraying coating treatment liquids or quenched by coolants, understanding the
    dynamics of liquid coating on fiber net becomes essential. For example, spray coating mesh technology is used in oil-water separation that
    is a worldwide problem due to the increasing emission of industrial oily waste water and the frequent oil spill accidents.  Technology for
     construction of  super-hydrophobic attapulgite coated mesh for gravity driven oil water separation fabricated by a  spray-coating process
     is described by \cite{spray-coating}.  Spray coating is also widely utilized for coating process of solar cells.

 For all the problems above the structure of the domains can be represented by graphs and the coating process
 in the lubrication limit can be modeled as liquid thin-film dynamics. The graph domain structure was already used in the analysis of certain fluid flows (see e.\,g. \cite{Gus85,Shaf98}). Mean curvature flows on graphs were studied in \cite{MeanCurvatureGraph} and diffuse interface PDE models on graphs were
 analyzed in \cite{Diffuse2012, GraphAllen}. Graph theory has applications in many different areas of science like: in computer graphics, internet tomography, quantum computing), physics (e.g., Anderson localization, photonic crystals, mesoscopic systems, waveguides), chemistry (aromatic molecules), and engineering (dynamical system, nanotechnology, microelectronics, fractal devices) (see for survey \cite{Kuch02,Leo10}).

 In the present paper, a coupled system of thin-film equations (shortly TFEs) with Kirchhoff's boundary conditions at the inner nodes and Neumann
boundary condition at the boundary nodes is used to describe viscous liquid coating of a fiber net. This model
was obtained as lubrication approximation of the Navier-Stokes system for incompressible flows.
The graphs can be interpreted as narrow  grooves on a solid surface in which extends a viscous fluid.
Our study allows to extend the previously obtained results (see \cite{BF90, nonnegative, multi, blow}) to the case of surfaces
with more complex geometry. To the best of our knowledge, this result is new and no other authors studied
TFEs on graphs previously.

Let us briefly describe the contents of the article.
In the next section we present graph structure of mesh domain, some definitions and auxiliary statements. In section 3, for the non-linearity power $n \geqslant 1$, we prove existence of non-negative weak solutions for TFEs on graph domains. The last section 4 is devoted to the proof of convergence toward a constant steady state. This section also includes numerical simulations of convergence to uniform coating for some different configurations of graphs.

\section{Notations and definitions}
\subsection{Graph structure of mesh domains}

Let $G=(V,E)$ be a metric graph with vertex set $V = \{a_i \}_{i=1}^m$, and the edge set
$E = \{e_j \}_{j=1}^l$ with $|e_j| = \ell_j$ and $e_j$ has the cross-sectional area $d_j > 0$.
Further, for simplicity, we will assume that $d_j = 1$ (one can introduce different weights $d_j$ to the edges to model more general geometry). Let $h(x,t)$ be a function defined
on $G \times \mathbb{R}^+$, $h_j(x,t)$ be its parameterization realization on $e_j \times \mathbb{R}^+$.
If $h_j(x,t)$, $j \in \{1,2,..,l\}$, satisfy the partial differential equation
\begin{equation}\label{a-1}
h_{j,\tau}  +  \bigl( h_{j}^n h_{j,xxx} \bigr)_x = 0 , \ \  x \in e_j:=(\alpha_j,\beta_j),
\end{equation}
where $\beta_j - \alpha_j = \ell_j > 0$, then $h(x,\tau)$ is called satisfying the TFE on $E$.
For a function $h(x,\tau)$ satisfied the TFE, we can define its normalized realization on $e_j$ by
$$
u_j(s,t) = h_j(\alpha_j + s\ell_j, \tau), \ \ s \in (0,1), \text{ and } \tau = \ell_j^4 t.
$$
Then we have
$$
u_{j,s} (s,t) = \ell_j  h_{j,x} (\alpha_j + s\ell_j,\tau),  \ \
\bigl( h_{j}^n h_{j,xxx} \bigr)_x = \ell_j^{-4} \bigl( u_{j}^n u_{j,sss} \bigr)_s,
$$
$$
u_{j,t} (s,t) = \ell_j^{-4}  h_{j,\tau} (\alpha_j + s\ell_j,\tau).
$$
So we can assume that $u_j(s,t)$ satisfies the TFE
\begin{equation}\label{a-1-0}
u_{j,t}  +  \bigl( u_{j}^n u_{j,sss} \bigr)_s = 0 , \ \  s \in (0,1), \ j =\overline{1,l}.
\end{equation}
The function $u_j(s,t)$ is called the normalized realization of $h(x,\tau)$. In the sequel, we always use
the normalized realization of a function.

At the interior node $a \in V_{int}:=V \setminus \partial G$ we assume that
\begin{equation}\label{inter}
\tfrac{\partial^k u_j(1,t)}{\partial s^k} = \tfrac{\partial^k u_i(0,t)}{\partial s^k}
\ \ \forall\,  j \in J^+(a), \ i \in J^-(a), \ k =0,2,
\end{equation}
\begin{equation}\label{ext-3}
\sum \limits_{j \in J^+(a)}  { \tfrac{\partial^k u_j(1,t)}{\partial s^k}   } -
\sum \limits_{j \in J^-(a)}  { \tfrac{\partial^k u_j(0,t)}{\partial s^k}  } =0 ,\ k = 1,3.
\end{equation}
Here (\ref{inter}) mean the nodal continuity of $u$ and its second derivatives or Kirchhoff's rules,
and  (\ref{ext-3}) are the flow continuous conditions.
At the boundary node $a \in \partial G$ we assume that
\begin{equation}\label{ext-1}
\tfrac{\partial^k u_j(1,t)}{\partial s^k} =   \tfrac{\partial^k u_i(0,t)}{\partial s^k} =0
\ \ \forall\,  j \in J^+(a), \ i \in J^-(a), \ k =1,3,
\end{equation}
where (\ref{ext-1}) are no-flux conditions. Thus, the corresponding closed loop system is
\begin{equation}\label{system}
\left \{\begin{gathered}
u_{j,t}  +   \bigl( u_{j}^n u_{j,sss} \bigr)_s = 0 , \  s \in (0,1),\ j =\overline{1,l},\\
\tfrac{\partial^k u_j(1,t)}{\partial s^k} = \tfrac{\partial^k u_i(0,t)}{\partial s^k}
\  \forall\,  j \in J^+(a), \ i \in J^-(a),
\ a \in V_{int}, \ k = 0,2,\\
\sum \limits_{j \in J^+(a)}  { \tfrac{\partial^k u_j(1,t)}{\partial s^k}   } -
\sum \limits_{j \in J^-(a)}  { \tfrac{\partial^k u_j(0,t)}{\partial s^k}  } =0 ,\ a \in V_{int}, \ k = 1,3, \\
\tfrac{\partial^k u_j(1,t)}{\partial s^k} =   \tfrac{\partial^k u_i(0,t)}{\partial s^k} =0
\ \ \forall\,  j \in J^+(a), \ i \in J^-(a), \ a \in \partial G,\ k =1,3,\\
u_j(s,0) = u_{0j}(s), \ \  s \in (0,1),\ j =\overline{1,l}.
\end{gathered} \right.
\end{equation}

\subsection{Functional spaces and definitions}

Define the function spaces $L^2(E)$ and $H^k(E)$ by
$$
L^2(E) = \{ f(x) : f_j(s) \in L^2(\alpha_j,\beta_j)\},
$$
$$
H^k(E) = \{ f(x) \in L^2(E) : f_j(s) \in H^k(\alpha_j,\beta_j)\},
$$
and the scalar product
$$
(u(x),v(x)) := \sum \limits_{j = 1}^l {\int \limits_0^1 {u_{j} (s ) v_j(s ) \,ds } }
$$
for arbitrary $u(x) = (u_1(x),..,u_l(x))$, $v(x) = (v_1(x),..,v_l(x)) \in L^2(E)$.

\begin{definition}
For node $a \in V$, let $J^+(a )$ denote the index set of the \textbf{incoming edges} to $a $ and $J^-(a )$
denote the index set of the \textbf{outgoing edges} from $a $, and
$$
\mathop {\lim} \limits_{s \to 1} u_j(s) = u_j(1 ) \text{ if } j \in  J^+(a ),\ \
\mathop {\lim} \limits_{s \to 0} u_j(s ) = u_j(0 ) \text{ if } j \in  J^-(a ),
$$
where $u_j$ is the normalized realization of $u(x)$ on $e_j$.
\end{definition}

\begin{definition}
A function $u(x)$ defined on $G$ is said to be the \textbf{incoming continuous} at
$ a \in V$ if $u(x)$ is continuous on $E$ and has limits at two endpoints of each edge in $E$, moreover
it satisfies
$$
u_j(1) = u(a) \ \ \forall\,  j \in J^+(a),
$$
where $u_j$ is the normalized realization of $u(x)$ on $e_j$. It is said to be the \textbf{outgoing continuous}
at $a$ if $u(x)$ is continuous on $E$ and has limits at two
endpoints of each edge in $E$, and
$$
u_i(0) = u(a) \ \ \forall\,  i \in J^-(a).
$$
For a multiple node $a$, $u(x)$ is said to be continuous at $a$ if $\mathop {\lim} \limits_{x \to a} u (x) = u(a)$
or equivalently
$$
 u(a) = u_j(1) = u_i(0) \ \ \forall\,  j \in J^+(a), \ i \in J^-(a).
$$
A function $u$ defined on $G$ is said to be a \textbf{continuous function} if it is continuous on $E$, and
continuous at each interior vertex $a \in V_{int}$, and at each boundary vertex $a_i \in \partial G$, it holds that
$$
\mathop {\lim} \limits_{s \to 1} u(s) = u_j(1 ) \text{ if } j \in  J^+(a_i),\ \
\mathop {\lim} \limits_{s \to 0} u(s ) = u_k(0 ) \text{ if } k \in  J^-(a_i ).
$$
One denotes the set of all continuous function on $G$ by $C(G)$.
\end{definition}

For more details about definitions in graph theory, see e.\,g. \cite{Xu10}.

\section{Main result}

Let us denote by
$$
G_0(z) : = \int \limits_A^z {\int \limits_A^v { \tfrac{dy dv }{|y|^n } } }   , \ \  A > 0 .
$$

\begin{theorem}\label{Th1}
Assume that $n \geqslant 1$ and
$$
0 \leqslant u_0(s) \in H^1(E), \ \  \sum \limits_{j = 1}^l { \int \limits_0^1 { G_0(u_{0j} (s ))  \,ds } } < + \infty.
$$
Let $G$ be a connected, simple, plane graph with $\partial G \neq \varnothing$.
Then there exists a nonnegative solution $u(s,t) =(u_1(s,t),..,u_l(s,t)) \in L^{\infty}(0,T;H^1(E))\cap L^{2}(0,T;H^2(E)) $ satisfying
$$
u_t \in L^{2}(0,T;(H^1(E))^*) , \  \sum \limits_{j = 1}^l {   \int \limits_0^T {
 \int \limits_0^1 {u_{j}^n (s,t) u_{j,sss}^2(s,t) \,ds } dt} } < \infty,
$$
$$
\sum \limits_{j = 1}^l { \int \limits_0^1 {u_j (s,t) \,ds } } = \sum \limits_{j = 1}^l { \int \limits_0^1 {u_{0j} (s) \,ds } } \qquad \mbox{(mass conservation)},
$$
and (\ref{system}) in the
following sense:
$$
\int \limits_0^T { < u_t, \psi >_{(H^1)^*, H^1} dt}  - \sum \limits_{j = 1}^l {   \int \limits_0^T {  \int \limits_0^1 {u^n_j(s,t)
u_{j,sss} (s,t) \psi_{j,s}(s,t) \,ds } dt}  } = 0
$$
for all $\psi (s,t)=(\psi_1(s,t),..,\psi_l(s,t)) \in L^{2}(0,T;H^1(E))$ and $T > 0$ such that
$$
\psi_j(1,t)  = \psi_i(0,t)  \  \forall\,  j \in J^+(a), \ i \in J^-(a),
\ a \in V.
$$
\end{theorem}

\subsection{Proof of Theorem~\ref{Th1}}

\subsubsection{Approximation solutions}
We write the approximation of the problem (\ref{system}) in the following form
\begin{equation}\label{system-2}
\left \{\begin{gathered}
u_{j,t}  -   \bigl( f_{\varepsilon}(u_{j}) w_{j,s} \bigr)_s = 0 , \  s \in (0,1),\ j =\overline{1,l},\\
w_j = -   u_{j,ss}, \ j =\overline{1,l},\\
u_j(1,t) = u_i(0,t)
\  \forall\,  j \in J^+(a), \ i \in J^-(a),
\ a \in V_{int},\\
w_j(1,t) = w_i(0,t) \  \forall\,  j \in J^+(a), \ i \in J^-(a),
\ a \in V_{int},\\
\sum \limits_{j \in J^+(a)}  {  u_{j,s}(1,t)   } - \sum \limits_{j \in J^-(a)}  {  u_{j,s}(0,t)  } =0,\ a \in V_{int},\\
\sum \limits_{j \in J^+(a )}  { w_{j,s}(1,t) }
- \sum \limits_{j \in J^-(a )}  { w_{j,s}(0,t) } = 0, \ a \in V_{int},\\
u_{j,s}(1,t)  =   u_{i,s}(0,t) = 0
\ \forall\,  j \in J^+(a), \ i \in J^-(a),
\ a \in \partial G,   \\
w_{j,s}(1,t)  =   w_{i,s}(0,t) = 0 \ \forall\,  j \in J^+(a), \ i \in J^-(a), \ a \in \partial G,\\
u_j(s,0) = u_{0j}^\varepsilon(s) \geqslant u_{0j}(s) + \varepsilon^{\theta}, \ \  s \in (0,1),\ j =\overline{1,l},
\end{gathered} \right.
\end{equation}
where $f_{\varepsilon}(z):= |z|^n + \varepsilon$, $\theta \in (0,\frac{1}{2})$.
To prove the local in time existence, we apply the Galerkin method. Let $\{ \phi_{ik} \}_{i,k=1}^{l,N}$
 be the eigenfunctions of the Laplace operator
$$
- \phi''_{ik}(s) = \lambda_{ik}   \phi_{ik} (s), \ \  s \in (0,1),\ i =\overline{1,N},\ k =\overline{1,l},
$$
with the continuity conditions
$$
\phi_{ik}(1)  =  \phi_{ij}(0) \  \forall\,  k \in J^+(a), \ j \in J^-(a), \ a \in V_{int},
$$
$$
\sum \limits_{k \in J^+(a)}  {  \phi' _{ik}(1) }
- \sum \limits_{k \in J^-(a)}  {  \phi' _{ik}(0) } = 0, \ a \in V_{int},
$$
$$
\phi'_{ik}(1)  =   \phi'_{ij}(0)=0 \  \forall\,  k \in J^+(a), \ j \in J^-(a),\  a \in \partial G.
$$
The eigenfunctions $\phi_{ik}$ are orthogonal in the $H^1(0,1)$ and orthonormal in the $L^2(0,1)$ scalar product,
i.\,e.
$$
 \int \limits_0^1 {\phi_{jk}(s) \phi_{ik}(s) \,ds } = 0 \text{ if } j \neq i,
\text{ and }  = 1 \text{ if } j =i;
$$
$$
 \int \limits_0^1 {\phi'_{jk}(s) \phi'_{ik}(s) \,ds }  = 0 \text{ if } j \neq i,
\text{ and }  = \lambda_{jk}   \text{ if } j =i.
$$
For more details about Sturm-Liouville theory on graphs, see e.\,g. \cite{Bel88,Bel8898,Pok05}.
Now, we consider the following Galerkin
ansatz
$$
u_j^{N,\varepsilon} (s,t) = \sum \limits_{i =1}^{N}  {  c_{ij}(t) \phi_{ij}(s) }, \ \
w_j^{N,\varepsilon} (s,t) = \sum \limits_{i =1}^{N}  {  d_{ij}(t) \phi_{ij}(s) }.
$$
Plugging this ansatz into $(\ref{system})_{1,2}$, multiplying by $\phi_{ij}(s)$, we obtain
$$
d_{ij}(t) =   \lambda_{ij} c_{ij}(t),
$$
\begin{equation}\label{c-1}
c'_{ij}(t) = -   \sum \limits_{k =1}^{N} { \lambda_{kj} c_{kj}(t)
\int \limits_0^1 {f_{\varepsilon}\Bigl( \sum \limits_{k =1}^{N}  {  c_{kj}(t) \phi_{kj}(s) } \Bigr)  \phi'_{kj}(s) \phi'_{ij}(s) \,ds } },
\end{equation}
\begin{equation}\label{c-2}
c_{ij}(0) = \int \limits_0^1 {u_{0j}^\varepsilon(s) \phi_{ij}(s) \,ds },
\end{equation}
which have to hold for $i = \overline{1,N}$, $j = \overline{1,l}$. Since the right-hand side of (\ref{c-1}) is Lipschitz continuous
on $c_{ij}$. Thus by the Picard-Lindel\"{o}f and Cauchy theorems a unique global in time solution of (\ref{c-1})--(\ref{c-2}) exists.

Global solvability for arbitrary but fixed $T > 0$ can be proved by using a priori estimates (uniformly in $N$ and $\varepsilon$) which
will be obtained in the following subsection.

\begin{example}
Let $G$ be a planar graph such that the directed edges are defined by
$$
e_1 =(a_1,a_5), \ e_2 =(a_2,a_6), \ e_3 =(a_3,a_7), \ e_4 =(a_4,a_8),
$$
$$
e_5 =(a_5,a_6), \ e_6 =(a_6,a_8), \ e_7 =(a_7,a_8), \ e_8 =(a_5,a_7),
$$
the boundary of $G$ is $\partial G = \{a_1,a_2, a_3, a_4 \}$ (see Figure 1).

\begin{figure}[t]\label{figure0}
\label{fig31}
\begin{center}
\includegraphics[height= 6cm]{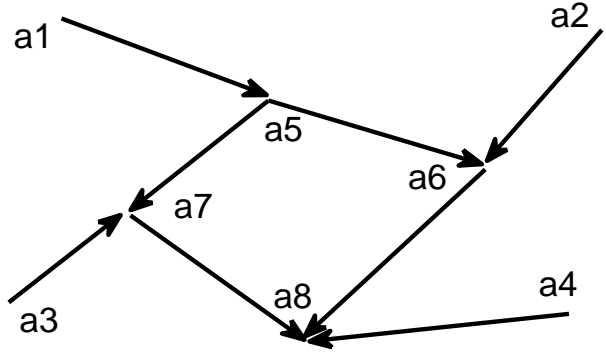}
\end{center}
\caption{An example of the graph with nonempty set of boundary nodes.}
\end{figure}

On the edges $e_k$, $k=1,2,3,4$, we have the problems
$$
- \phi''_{ik}(s) = \lambda_{ik} \phi_{ik}(s), \ \phi'_{ik}(0) = \phi'_{ik}(1) = 0, \ k = \overline{1,4}.
$$
The solutions are
$$
\phi_{ik}(s) = \sqrt{2} \cos( \sqrt{\lambda_{ik}} s), \ \lambda_{ik} = ( \pi i)^2, \ k = \overline{1,4}, \ i=\overline{1,N}.
$$
On the edges $e_k$, $k=5,6,7,8$, we get the problems
$$
- \phi''_{i5}(s) = \lambda_{i5} \phi_{i5}(s), \  \phi_{i1}(1)=\phi_{i5}(0) = \phi_{i8}(0),\ \phi'_{i5}(0) = -\phi'_{i8}(0);
$$
$$
- \phi''_{i6}(s) = \lambda_{i6} \phi_{i6}(s), \  \phi_{i2}(1)=\phi_{i6}(0) = \phi_{i5}(1),\ \phi'_{i6}(0) =  \phi'_{i5}(1);
$$
$$
- \phi''_{i7}(s) = \lambda_{i7} \phi_{i7}(s), \  
\phi'_{i7}(1) = - \phi'_{i6}(1);
$$
$$
- \phi''_{i8}(s) = \lambda_{i8} \phi_{i8}(s), \  \phi_{i3}(1)=\phi_{i8}(1) = \phi_{i7}(0),\ \phi'_{i8}(1) =  \phi'_{i7}(0);
$$
The corresponding solutions are
$$
\phi_{i5}(s) = (-1)^i \sqrt{2} \cos( \sqrt{\lambda_{i5}} s) 
, \ \lambda_{i5} = (2 \pi i)^2, \ i=\overline{1,N},
$$
$$
\phi_{i6}(s) = (-1)^i \sqrt{2} \cos( \sqrt{\lambda_{i6}} s) 
, \ \lambda_{i6} = (\pi i)^2, \ i=\overline{1,N},
$$
$$
\phi_{i7}(s) = (-1)^i \sqrt{2} \cos( \sqrt{\lambda_{i7}} s) 
, \ \lambda_{i7} = (\pi i)^2, \ i=\overline{1,N},
$$
$$
\phi_{i8}(s) = (-1)^i \sqrt{2} \cos( \sqrt{\lambda_{i8}} s) 
, \ \lambda_{i8} = (2 \pi i)^2, \ i=\overline{1,N}.
$$
%
\end{example}

\subsubsection{A priori estimates}

Next, for brevity, we denote by $u_j := u_j^{N,\varepsilon}$.
Integrating (\ref{system-2})$_1$ on $s$ and summing on $j$, we find that
\begin{multline*}
\tfrac{d}{dt} \sum \limits_{j = 1}^l { \int \limits_0^1 {u_j (s,t) \,ds } } =
\sum \limits_{j = 1}^l { f_\varepsilon( u_j (1,t)) w_{j,s}(1,t) }  -
 \sum \limits_{j = 1}^l {  f_\varepsilon( u_j (0,t)) w_{j,s}(0,t) }  =  \\
\sum \limits_{k = 1}^m \Bigl[  \sum \limits_{j \in J^+(a_k)}  {  f_\varepsilon( u_j (1,t)) w_{j,s }(1,t) }
- \sum \limits_{j \in J^-(a_k)}  { f_\varepsilon(  u_j (0,t)) w_{j,s }(0,t) }   \Bigr]  \mathop {=} \limits^{(\ref{system-2})_3}\\
\sum \limits_{k = 1}^m  f_\varepsilon( u(a_k ,t)) \Bigl[  \sum \limits_{j \in J^+(a_k)}  {   w_{j,s }(1,t) }
- \sum \limits_{j \in J^-(a_k)}  {     w_{j,s }(0,t) }   \Bigr] \mathop {=} \limits^{(\ref{system-2})_6} 0,
\end{multline*}
whence  we get
\begin{equation}\label{mass}
\sum \limits_{j = 1}^l { \int \limits_0^1 {u_j (s,t) \,ds } } = \sum \limits_{j = 1}^l { \int \limits_0^1 {u_{0j}^{\varepsilon} (s) \,ds } }.
\end{equation}
The energy function is defined by
$$
\mathcal{E}_{\varepsilon}(t) := \tfrac{1}{2} \sum \limits_{j = 1}^l { \int \limits_0^1 {u_{j,s}^2 (s,t) \,ds } }.
$$
Thus
\begin{multline*}
\tfrac{d \mathcal{E}_{\varepsilon}(t)}{dt} = \sum \limits_{j = 1}^l { \int \limits_0^1 {u_{j,s}(s,t) u_{j,st}(s,t) \,ds } } =
  \sum \limits_{j = 1}^l { \int \limits_0^1 {w_{j }(s,t) u_{j,t}(s,t) \,ds } }  +\\
\sum \limits_{k = 1}^m \Bigl[  \sum \limits_{j \in J^+(a_k)}  {  u_{j,s}(1,t) u_{j,t}(1,t) }
- \sum \limits_{j \in J^-(a_k)}  {  u_{j,s}(0,t) u_{j,t}(0,t) }   \Bigr]   =\\
\sum \limits_{k = 1}^m \Bigl[  \sum \limits_{j \in J^+(a_k)}  {  f_\varepsilon( u_{j}  (1,t)) w_{j }(1,t) w_{j,s }(1,t) }
- \sum \limits_{j \in J^-(a_k)}  {  f_\varepsilon( u_{j}  (0,t)) w_{j }(0,t) w_{j,s }(0,t) }   \Bigr]  +\\
\sum \limits_{k = 1}^m \Bigl[  \sum \limits_{j \in J^+(a_k)}  {  u_{j,s}(1,t) u_{j,t}(1,t) }
- \sum \limits_{j \in J^-(a_k)}  {  u_{j,s}(0,t) u_{j,t}(0,t) }   \Bigr] -\\
\sum \limits_{j = 1}^l {   \int \limits_0^1 { f_\varepsilon(u_{j}  (s,t)) w_{j,s }^2(s,t) \,ds } } \mathop {=} \limits^{(\ref{system-2})_{3,4}}\\
\sum \limits_{k = 1}^m f_\varepsilon( u (a_k,t)) w(a_k,t)\Bigl[  \sum \limits_{j \in J^+(a_k)}  {   w_{j,s }(1,t) }
- \sum \limits_{j \in J^-(a_k)}  {   w_{j,s }(0,t) }   \Bigr]  +\\
\sum \limits_{k = 1}^m   u_{t}(a_k,t) \Bigl[  \sum \limits_{j \in J^+(a_k)}  {  u_{j,s}(1,t)   }
- \sum \limits_{j \in J^-(a_k)}  {  u_{j,s}(0,t)  }   \Bigr] - \\
\sum \limits_{j = 1}^l {   \int \limits_0^1 { f_\varepsilon(u_{j} (s,t)) w_{j,s }^2(s,t) \,ds } },
\end{multline*}
whence, due to (\ref{system-2})$_5$ and (\ref{system-2})$_6$, we obtain
\begin{equation}\label{energy}
 \mathcal{E}_\varepsilon(t)  + \sum \limits_{j = 1}^l {   \int \limits_0^t {
 \int \limits_0^1 {f_\varepsilon(u_{j}  (s,t)) w_{j,s }^2(s,t) \,ds } dt} } =  \mathcal{E}_\varepsilon(0),
\end{equation}
hence $ \mathcal{E}_\varepsilon(t) \leqslant  \mathcal{E}_\varepsilon(0) $. This means that the energy of the closed loop system
(\ref{system-2}) is dissipative.


The entropy function is defined by
$$
G_\varepsilon(u) :=     \int \limits_A^u {\int \limits_A^v { \tfrac{dy dv }{f_\varepsilon(y)}   } },
\ \ G''_\varepsilon(u) = \tfrac{1}{f_\varepsilon(u)}   \geqslant 0, \  A > 0.
$$
Thus
\begin{multline*}
\tfrac{d}{dt}  \sum \limits_{j = 1}^l { \int \limits_0^1 { G_\varepsilon(u_{j} (s,t))  \,ds } } = \\
 \sum \limits_{j = 1}^l { \int \limits_0^1 {G'_\varepsilon(u_{j} (s,t)) u_{j,t}(s,t) \,ds } } =
- \sum \limits_{j = 1}^l { \int \limits_0^1 { u_{j,s}(s,t)w_{j,s}(s,t) \,ds } }  +\\
   \sum \limits_{k = 1}^m \Bigl[  \sum \limits_{j \in J^+(a_k)}  {  G'_\varepsilon( u_{j}(1,t)) f_\varepsilon(u_{j}(1,t))  w_{j,s }(1,t) }
\\  - \sum \limits_{j \in J^-(a_k)}  {  G'_\varepsilon( u_{j}(0,t)) f_\varepsilon(u_{j}(0,t)) w_{j,s }(0,t) }   \Bigr]   =
- \sum \limits_{j = 1}^l { \int \limits_0^1 {  w_{j }^2(s,t) \,ds } } -\\
\sum \limits_{k = 1}^m \Bigl[  \sum \limits_{j \in J^+(a_k)}  {    u_{j,s}(1,t) w_{j }(1,t) }
- \sum \limits_{j \in J^-(a_k)}  {  u_{j,s}(0,t) w_{j }(0,t) }   \Bigr] +\\
\sum \limits_{k = 1}^m \Bigl[  \sum \limits_{j \in J^+(a_k)}  {  G'_\varepsilon( u_{j}(1,t)) f_\varepsilon(u_{j}(1,t))  w_{j,s }(1,t) }
\\ - \sum \limits_{j \in J^-(a_k)}  {  G'_\varepsilon( u_{j}(0,t)) f_\varepsilon(u_{j}(0,t)) w_{j,s }(0,t) }   \Bigr]  \mathop {=} \limits^{(\ref{system-2})_{3,4}}
- \sum \limits_{j = 1}^l { \int \limits_0^1 {  w_{j }^2(s,t) \,ds } } -\\
\sum \limits_{k = 1}^m w(a_k,t) \Bigl[  \sum \limits_{j \in J^+(a_k)}  {    u_{j,s}(1,t)  }
- \sum \limits_{j \in J^-(a_k)}  {   u_{j,s}(0,t)   }   \Bigr]+\\
  \sum \limits_{k = 1}^m G'_\varepsilon(u (a_k,t)) f_\varepsilon(u (a_k,t))\Bigl[  \sum \limits_{j \in J^+(a_k)}  {    w_{j,s }(1,t) }
- \sum \limits_{j \in J^-(a_k)}  {   w_{j,s }(0,t) }   \Bigr],
\end{multline*}
whence, due to (\ref{system-2})$_5$ and (\ref{system-2})$_6$, we obtain
\begin{equation}\label{entropy}
\sum \limits_{j = 1}^l { \int \limits_0^1 { G_\varepsilon(u_{j} (s,t))  \,ds } } + \sum \limits_{j = 1}^l {   \int \limits_0^t {
 \int \limits_0^1 { w_{j }^2(s,t) \,ds } dt} } =  \sum \limits_{j = 1}^l { \int \limits_0^1 { G_\varepsilon(u_{0j}^{\varepsilon} (s ))  \,ds } }.
\end{equation}
This means that the entropy of the closed loop system
(\ref{system-2}) decays.

As a result, in view of $u_0 \in H^1(E)$, from (\ref{mass}), (\ref{energy}) and (\ref{entropy}) we obtain that
$$
\{u^{N,\varepsilon}\}_{N\geqslant 1, \varepsilon >0} \text{ is uniformly bounded in } L^\infty(0,T; H^1(E)),
$$
$$
\{w^{N,\varepsilon} \}_{N\geqslant 1, \varepsilon >0} \text{ is uniformly bounded in } L^2(Q_T),
$$
$$
\{u_t^{N,\varepsilon} \}_{N\geqslant 1, \varepsilon >0} \text{ is uniformly bounded in } L^2(0,T; (H^1(E))^*),
$$
$$
\{G_\varepsilon (u^{N,\varepsilon})\}_{N\geqslant 1, \varepsilon >0} \text{ is uniformly bounded in } L^\infty(0,T; L^1(E)),
$$
$$
\{(u^{N,\varepsilon})^{\frac{n}{2}} w^{N,\varepsilon}_s\}_{N\geqslant 1, \varepsilon >0} \text{ is uniformly bounded in } L^2(Q_T).
$$
Following \cite{BF90}, we can let $N \to +\infty,\, \varepsilon \to 0$, and prove
nonnegativity of $u(s,t)$ for $n \geqslant 1$. As a result, the proof of Theorem 1 is complete.

\begin{figure}[t]\label{fig32}
\begin{center}
\hskip -1 cm
\includegraphics[height= 6cm]{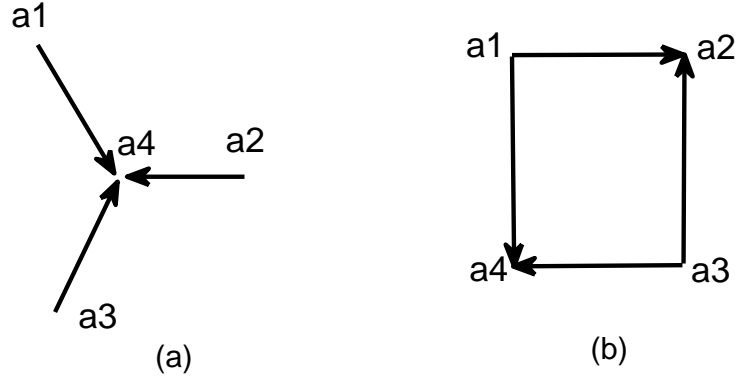}
\end{center}
\caption{The graph structure of the mesh domain that is used for numerical simulations of convergence to uniform coating for 3-edges case with non-empty set of boundary nodes (a) and for 4-edges case with an empty set of boundary nodes (b).}
\end{figure}

\begin{figure}[t]
\label{fig33}
\centering
\includegraphics[scale=.65]{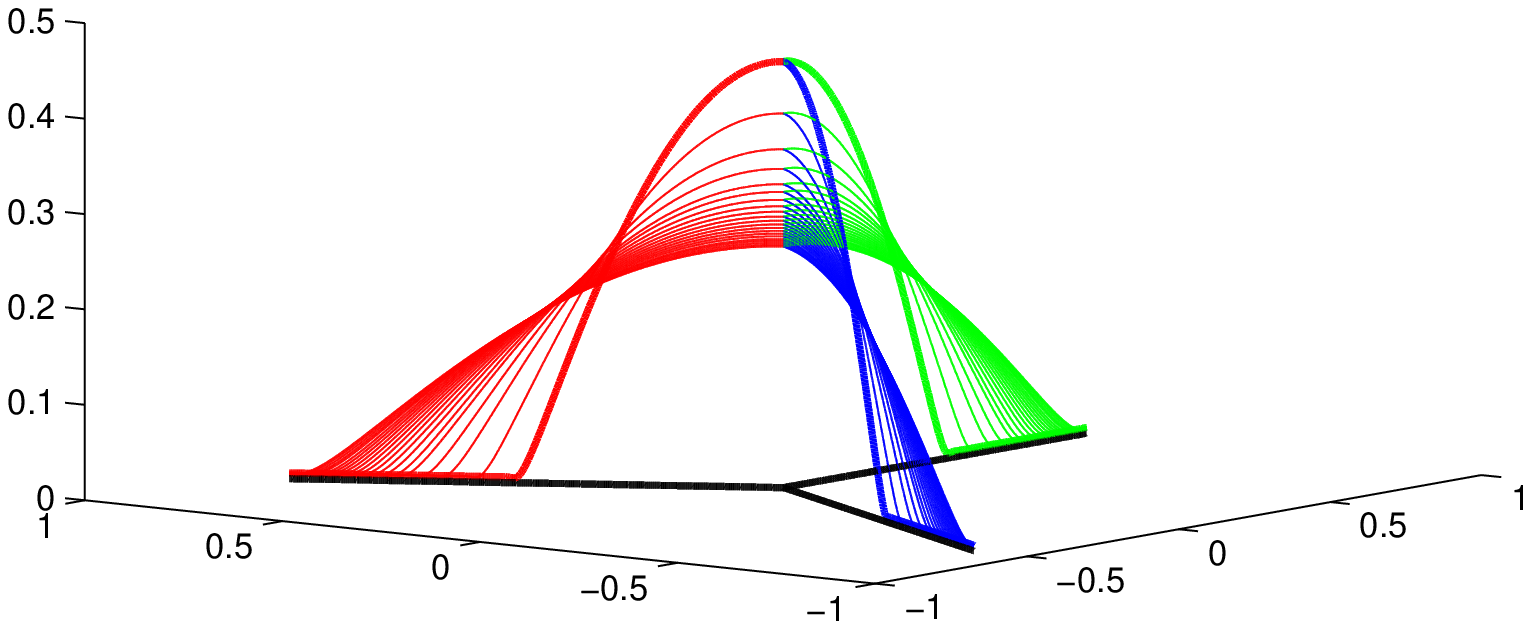}
\vskip 0.2cm
\includegraphics[scale=.65]{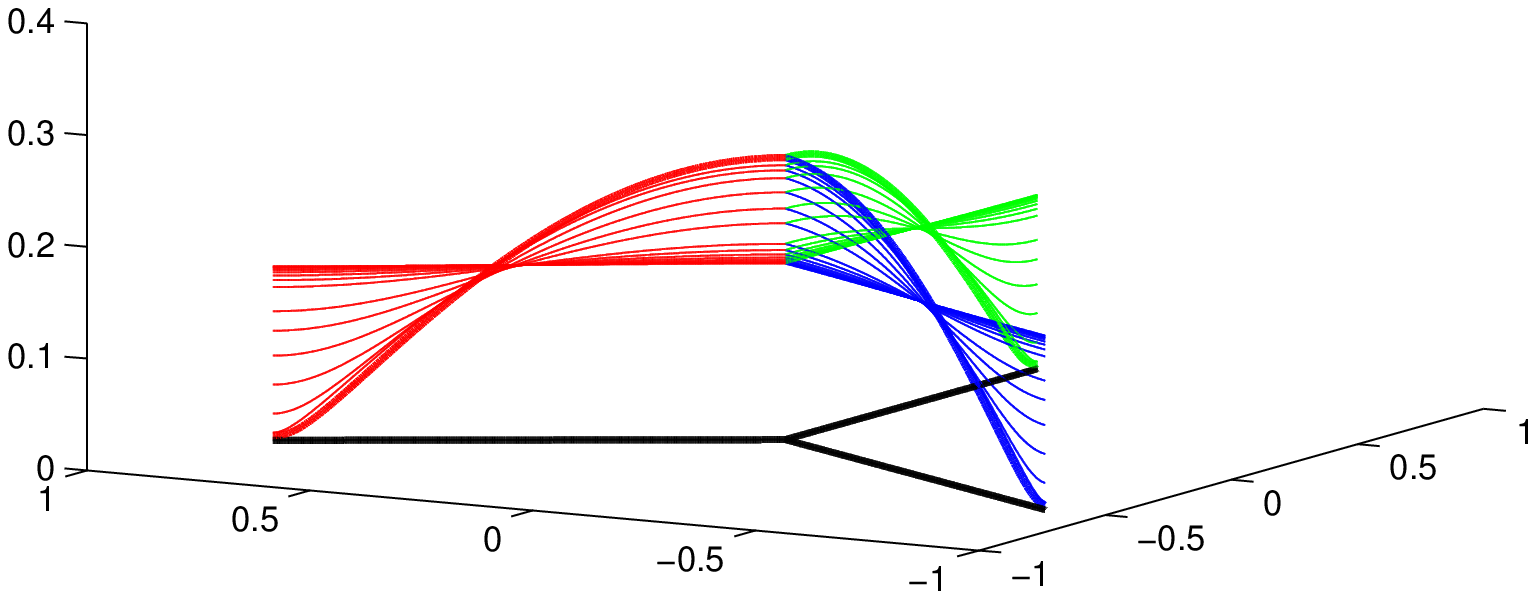}
\caption{Snapshots of numerical time evolution for coating flow on 3-edge graph domain (symmetric initial values). Local time dynamics (on the top) and long time dynamics (on the bottom). }
\end{figure}



\section{Convergence to steady state and numerical simulations}

By the Cauchy inequality we have
$$
  \int \limits_{s_0}^{s} { u_{j} (s,t)  w_{j,s } (s,t) \,ds } \leqslant
\Bigl( \int \limits_0^1 {f_\varepsilon(u_{j}  (s,t)) w_{j,s }^2(s,t) \,ds }\Bigr)^{\frac{1}{2}}
\Bigl( \int \limits_0^1 {\tfrac{u_{j}^2 (s,t)}{f_\varepsilon(u_{j}(s,t))} \,ds }\Bigr)^{\frac{1}{2}}
$$
for all $s_0,\,s \in (0,1)$, and $j = \overline{1,l} $. On the other hand,
\begin{multline*}
 \int \limits_{s_0}^{s} { u_{j} (s,t)  w_{j,s } (s,t) \,ds } = u_{j} (s,t)  w_{j} (s,t) +
\tfrac{1}{2} u_{j,s}^2(s,t) -    \\
u_{j} (s_0,t)  w_{j} (s_0,t) -
\tfrac{1}{2} u_{j,s}^2(s_0,t), \text{ where } w_j = - u_{j,ss}.
  \end{multline*}
From here we get
\begin{multline*}
u_{j} (s,t)  w_{j} (s,t) + \tfrac{1}{2} u_{j,s}^2(s,t) - f_j(s_0,t) \leqslant    \\
 \Bigl( \int \limits_0^1 {f_\varepsilon(u_{j}  (s,t)) w_{j,s }^2(s,t) \,ds }\Bigr)^{\frac{1}{2}}
\Bigl( \int \limits_0^1 {\tfrac{u_{j}^2 (s,t)}{f_\varepsilon(u_{j}(s,t))} \,ds }\Bigr)^{\frac{1}{2}},
\end{multline*}
where $f_j(s_0,t) : = u_{j} (s_0,t)  w_{j} (s_0,t) + \tfrac{1}{2} u_{j,s}^2(s_0,t)$.
Integrating in $s$ over $(0,1)$, after summing on $j$, gives
\begin{multline*}
 \tfrac{3}{2} \sum \limits_{j = 1}^l { \int \limits_0^1 {u_{j,s}^2 (s,t) \,ds } } -
 \sum \limits_{k = 1}^m \Bigl[  \sum \limits_{j \in J^+(a_k)}  { \hspace{-0.5cm} u_{j }(1,t) u_{j,s}(1,t) }
- \sum \limits_{j \in J^-(a_k)}  { \hspace{-0.5cm}  u_{j }(0,t) u_{j,s}(0,t) }   \Bigr] -
\sum \limits_{j = 1}^l { f_j(s_0,t) } \leqslant \\
\sum \limits_{j = 1}^l  { \Bigl( \int \limits_0^1 {f_\varepsilon(u_{j}  (s,t)) w_{j,s }^2(s,t) \,ds }\Bigr)^{\frac{1}{2}}
\Bigl( \int \limits_0^1 {\tfrac{u_{j}^2 (s,t)}{f_\varepsilon(u_{j}(s,t))} \,ds }\Bigr)^{\frac{1}{2}} } \leqslant \\
\Bigl( \sum \limits_{j = 1}^l  { \int \limits_0^1 {f_\varepsilon(u_{j}  (s,t)) w_{j,s }^2(s,t) \,ds } }\Bigr)^{\frac{1}{2}}
\Bigl( \sum \limits_{j = 1}^l  { \int \limits_0^1 {\tfrac{u_{j}^2 (s,t)}{f_\varepsilon(u_{j}(s,t))} \,ds } }\Bigr)^{\frac{1}{2}}.
\end{multline*}
It follows from the boundary conditions in the problem (\ref{system}) that there exists $s_0 \in [0,1]$ such that $f_j(s_0,t) \leqslant 0$. So,
by (\ref{entropy}) we deduce that
\begin{multline}\label{w-1}
 9 \mathcal{E}_\varepsilon^2 (t)   \leqslant \Bigl( \sum \limits_{j = 1}^l  { \int \limits_0^1 {f_\varepsilon(u_{j}  (s,t)) w_{j,s }^2(s,t) \,ds } }\Bigr) \Bigl( \sum \limits_{j = 1}^l  { \int \limits_0^1 {\tfrac{u_{j}^2 (s,t)}{f_\varepsilon(u_{j}(s,t))} \,ds } }\Bigr) \leqslant \\
\Bigl( \sum \limits_{j = 1}^l  { \int \limits_0^1 {f_\varepsilon(u_{j}  (s,t)) w_{j,s }^2(s,t) \,ds } }\Bigr) \Bigl( \sum \limits_{j = 1}^l  { \int \limits_0^1 { |u_{j} (s,t)|^{2-n} \,ds } }\Bigr) \leqslant \\
C(u_{0 j}^{\varepsilon})\sum \limits_{j = 1}^l  { \int \limits_0^1 {f_\varepsilon(u_{j}  (s,t)) w_{j,s }^2(s,t) \,ds } }.
\end{multline}
From (\ref{energy}), due to (\ref{w-1}), we arrive at
$$
\frac{d}{dt}\mathcal{E}_\varepsilon  (t)  + \tfrac{9}{C(u_{0 j}^{\varepsilon})} \mathcal{E}_\varepsilon^2 (t) \leqslant 0,
$$
whence
$$
\mathcal{E}_\varepsilon  (t) \leqslant  \mathcal{E}_\varepsilon (0) \bigl( 1 +  \tfrac{9}{C(u_{0 j}^{\varepsilon})} \mathcal{E}_\varepsilon (0) t \bigr)^{-1}.
$$
Passing to the limit as $\varepsilon \to 0$, we obtain
\begin{equation}\label{w-2}
\mathcal{E}_0  (t) \leqslant  \mathcal{E}_0 (0) \bigl( 1 +  \tfrac{9}{C(u_{0 j})} \mathcal{E}_0 (0) t \bigr)^{-1} \to 0
\text{ as } t \to +\infty.
\end{equation}
As a result, $u_{j,s}(s,t) \to 0$ as $t \to +\infty$, by continuity $u (x,t)$ in each vertex, implies 
$u_{j }(s,t) \to K$ for all $j =\overline{1,l}$, where $K > 0$ is some constant. By the mass conservation
(\ref{mass}) with $\varepsilon = 0$, we find that
$K = \frac{1}{l}\sum \limits_{j = 1}^l { \int \limits_0^1 {u_{0j}  (s) \,ds } }$ .
Hence, we obtained the following result.
\begin{corollary}
For any $j = \overline{1,l}$
$$
u_j(s,t) \to \tfrac{1}{l}\sum \limits_{j = 1}^l { \int \limits_0^1 {u_{0j}  (s) \,ds } } \text{ as }
t \to + \infty.
$$
\end{corollary}

Two different types of graph domains, which were used in numerical simulations described below, are illustrated in Figure 2. For the case (a) we ran Matlab finite element numerical simulations for symmetric initial values (see Figure 3) and for non-symmetric initial values (see Figure 4). Neumann (no-flux) boundary conditions were used at $3$ boundary nodes $u_{1,x}=u_{2,x}=u_{3,x}= 0$,
$u_{1,xxx}= u_{2,xxx}= u_{3,xxx} = 0$ and Kirchhoff's boundary conditions were applied at the only inner node $u_{1,x} + u_{2,x}+ u_{3,x} = 0$,  $u_{1,xxx} + u_{2,xxx}+ u_{3,xxx} = 0$ with continuity conditions
$u_1= u_2 =u_3$, $u_{1,xx}= u_{2,xx}= u_{3,xx}$. On the top pictures (Figure 3, 4) bold lines are used to indicate initial data (for all edges initial data are given by droplet concentrated near the inner node). For local (short time dynamics) in both cases the initial droplets spread over their edges with the only difference that in case of symmetry all $3$ first derivatives at the inner node are equal to $0$. The last values of the numerical short time dynamics are used as initial values (see bold lines on the bottom pictures in Figure 3, 4) for long time dynamics time evolution snapshots. This long time numerics clearly illustrates convergence toward uniform coating in both (symmetric and non-symmetric) cases.

For the case (b) with an empty set of boundary nodes and non-symmetric initial values (see Figure 5). Kirchhoff's boundary conditions were applied at the $4$ inner nodes: $u_{1,x} + u_{2,x} =0$, $u_{2,x} + u_{3,x} = 0$, $u_{3,x} + u_{4,x} = 0$, $u_{4,x} + u_{1,x} = 0$, $u_{1,xxx} + u_{2,xxx} =0$,
$u_{2,xxx} + u_{3,xxx} = 0$, $u_{3,xxx} + u_{4,xxx} = 0$, $u_{4,xxx} + u_{1,xxx} = 0$ with corresponding continuity conditions $u_1=u_2, \, u_{1,xx}=u_{2,xx}$,  $u_2=u_3, \, u_{2,xx}=u_{3,xx}$, $u_3=u_4, \, u_{3,xx}=u_{4,xx}$, and $u_4=u_1, \, u_{4,xx}=u_{1,xx}$. On the top pictures (Figure 5) bold lines are used to indicate initial data (for edges $1$ and $3$ (blue and green) initial data are given by bigger droplets to compare to the edges $2$ and $4$ (yellow and red)). The long time numerical simulations show the difference in convergence toward uniform coating between the edges. On the edges $1$ and $3$
solutions approach the constant value from above and at the same time on the edges $2$ and $4$ solutions approach the constant value from below.

\vspace{5mm} {\it {\bf Acknowledgement.} This work was partially
supported by a grant from the Simons Foundation (\#275088 to Marina
Chugunova) and  by a grant from Ministry of Education and Science of Ukraine
(0118U003138 to Roman Taranets)}

\end{document}